\documentclass[12pt]{amsart}
\usepackage{amsfonts}
\usepackage{latexsym,amsmath,amsthm,amssymb,amsfonts}

\numberwithin{equation}{section}
\allowdisplaybreaks

\def\endproof{$\hfill\Box$\\}

\numberwithin{equation}{section}
\newtheorem{theorem}{Theorem}[section]
\newtheorem{lem}[theorem]{Lemma}
\newtheorem{thm}[theorem]{Theorem}

\newcounter{Cnumber}

\title[ ]
{\bf Local strong solution to General Landau-Lifshitz-Bloch equation }
\author[ ]
{Zonglin Jia}

\begin{document}
\maketitle

\begin{abstract}
In this paper, we bring in General Landau-Lifshitz-Bloch equation and prove that it admits a local strong solution.
\end{abstract}

\section{Introduction}
Landau-Lifshitz-Gilbert equation describes physical properties of micromagnetic at temperatures below the critical temperature. The equation is as follows:
\begin{equation}\label{1.1}
\frac{\partial m}{\partial t}=\lambda_1m\times H_{eff}-\lambda_2m\times(m\times H_{eff})
\end{equation}
where $"\times"$ denotes the vector cross product in $\mathbb{R}^3$ and $"H_{eff}"$ is effective field while $\lambda_1$ and $\lambda_2$ are real constants. From (\ref{1.1}) we can see that, if $M$ is a domain in $\mathbb{R}^d$, $m:M\times[0,T]\longrightarrow\mathbb{R}^3$ is a classical solution for (\ref{1.1}) and $m(x,0):=m_0(x)\in S^2$, then $m(x,t)\in S^2$(the sphere in $\mathbb{R}^3$). Indeed, one can multiply two sides of (\ref{1.1}) by $m$ to get that
\[\frac{d}{dt}|m|^2=0\]
So we always require that the solution to (\ref{1.1}) lies on $S^2$.

However, at high temperature, the model must be replaced by following Landau-Lifshitz-Bloch equation(LLB)
\begin{equation}\label{1.2}
\frac{\partial u}{\partial t}=\gamma u\times H_{eff}+L_1\frac{1}{|u|^2}(u\cdot H_{eff})u-L_2\frac{1}{|u|^2}u\times(u\times H_{eff})
\end{equation}
where $\gamma$, $L_1$, $L_2$ are real numbers and $\gamma>0$. $H_{eff}$ is given by
\[H_{eff}=\Delta u-\frac{1}{\chi||}\Big(1+\frac{3T}{5(T-T_c)}|u|^2\Big)u.\]
where $T>T_c>0$ and $\chi||>0$. Now we are not able to say that the classical solution $u(x,t)$ is always in $S^2$ because
\[\frac{1}{2}\frac{d}{dt}|u|^2=L_1(u\cdot H_{eff})\not=0\]

To the best of our knowledge the analysis of the LLB equation is an open problem at present. In \cite{L}, Le consider the case that $L_1=L_2=:\kappa_1>0$. At that time, he rewrites (\ref{1.2}) as
\[\frac{\partial u}{\partial t}=\kappa_1\Delta u+\gamma u\times\Delta u-\kappa_2(1+\mu|u|^2)u\]
with $\kappa_2:=\frac{\kappa_1}{\chi||}$ and $\mu:=\frac{3T}{5(T-T_c)}$ and assume that $\kappa_2$, $\gamma$, $\mu$ is positive. Le has proven that above equation with Neumann boundary value conditions has global weak solution(the "weak solution" here is different from ordinary one).\\

Inspired by Le, we bring in following equation
\begin{equation}\label{1.3}
\left\{
\begin{array}{llll}
\partial_t u=\kappa_1\Delta u+\gamma\nabla F(u)\times\Delta u-\kappa_2(1+\mu\cdot F(u))\nabla F(u)\\
\frac{\partial u}{\partial\nu}=0\\
u(\cdot ,0)=u_0
\end{array}
\right.
\end{equation}
where $\nu$ is outer normal direction of $\partial M$ and we call it Generalized Landau-Lifshitz-Bloch equation(GLLB). Using Galerkin approximation and introducing auxiliary functions
\[I(\lambda):=\sup\limits_{|z|\leqslant\lambda}|(HessF)(z)|\]
\[J(\lambda):=\sup\limits_{|z|\leqslant\lambda}|\nabla F(z)|\]
\[H(\lambda):=\sup\limits_{|z|\leqslant\lambda}|F(z)|\]
we get the upper bound of $L^{\infty}([0,T^*],W^{2,2}(M,\mathbb{R}^3))-$norm of approximation solutions for some $T^*$. Our main result is as follows:
\begin{thm}
Suppose that $F\in C^3(\mathbb{R}^3)$, $u_0\in W^{2,2}(M,\mathbb{R}^3)$ and $\frac{\partial u_0}{\partial\nu}=0$. $M$ is a regular bounded domain of $\mathbb{R}^d(d\leqslant3)$. then (\ref{1.3}) admits a strong solution $u\in
L^{\infty}([0,T^*],W^{2,2}(M,\mathbb{R}^3))$ for some $T^*>0$.
\end{thm}
\section{The Proofs of Theorem }
We appoint that $"||\cdot||_p"$ means $"||\cdot||_{L^p}"$. Before dealing with the main theorem, we need following useful lemma.
\begin{lem}
Suppose that $F:\mathbb{R}^m\longrightarrow\mathbb{R}^n$ is $C^1-$ smooth. Let
\[g(\lambda)=\max\limits_{|x|\leqslant\lambda}\{|F(x)|\},\]
then g is monotonously increasing and locally Lipschitz.
\end{lem}
\proof Obviously, $g$ is increasing.

Now taking any $\lambda<\mu$, we assume that
\[|F(x_{\lambda})|=g(\lambda)\]
and
\[|F(x_{\mu})|=g(\mu)\]
By the definition, the case when $|x_{\mu}|\leqslant\lambda$ is trival. So we only consider that $\lambda<|x_{\mu}|\leqslant\mu$. Let
\[x_0:=\lambda\frac{x_{\mu}}{|x_{\mu}|}\]
Since
\[|F(x_0)|\leqslant|F(x_{\lambda})|\]
we get
\[
\aligned
&g(\mu)-g(\lambda)\\
\leqslant&|F(x_{\mu})|-|F(x_0)|\\
\leqslant&|F(x_{\mu})-F(x_0)|\\
\leqslant&\max\limits_{|z|\leqslant\mu}\{|DF(z)|\}|x_{\mu}-x_0|\\
\leqslant&\max\limits_{|z|\leqslant\mu}\{|DF(z)|\}(\mu-\lambda)
\endaligned
\]
this complete the proof.\endproof\\

Now let us turn to theorem 1.1. From theorem 2.4.5 in chapter 2.4 of \cite{W}, there exists an orthonormal basis $\{e_i\}$ of $L^2(M)$ such that
\[-\Delta e_i=\lambda_i e_i\]
with
\[\frac{\partial e_i}{\partial\nu}=0\]
Let
\[u_n(t,x):=\sum\limits_{i=1}^nC_i^n(t)e_i(x)\]
where $C_i^n(t)$ is to be determined such that following equation holds
\begin{equation}\label{2.1}
\left\{
\begin{array}{llll}
\aligned
\int_M\partial_tu_n\cdot e_i\,dM=&\kappa_1\int_M\Delta u_n\cdot e_i\,dM+\gamma\int_M\nabla F(u_n)\times\Delta u_n\cdot e_i\,dM\\
&-\kappa_2\int_M(1+\mu\cdot F(u_n))\nabla F(u_n)\cdot e_i\,dM\\
u_n(0,\cdot)=u_{0n}
\endaligned
\end{array}
\right.
\end{equation}
where
\[u_{0n}:=\sum\limits_{i=1}^n\Big(\int_Mu_0\cdot e_i\,dM\Big)e_i.\]
Since (\ref{2.1}) is equivalent to an ODE in $\mathbb{R}^{3n}$ and the right-hand side of (\ref{2.1}) is locally Lipschitz respect to $C^n_i$, we can find the needed $C^n_i$ which exists for a short time. Using $C^n_i(t)$ to multiply two sides of (\ref{2.1}) and then summing $i$ from 1 to $n$, we get that
\[
\aligned
&\frac{1}{2}\frac{d}{dt}\int_M|u_n|^2\,dM+\kappa_1\int_M|\nabla u_n|^2\,dM\\
=&-\gamma\int_M(\nabla F(u_n)\times u_n)\Delta u_n\,dM-\kappa_2\int_M(1+\mu\cdot F(u_n))\nabla F(u_n)\cdot u_n\,dM\\
\leqslant&\gamma\int_M|\nabla F(u_n)|\cdot|u_n|\cdot|\Delta u_n|\,dM+\kappa_2\int_M(1+\mu|F(u_n)|)\cdot|\nabla F(u_n)|\cdot|u_n|\,dM\\
\leqslant&\gamma||\nabla F(u_n)||_{\infty}||u_n||_{\infty}\int_M|\Delta u_n|\,dM+\kappa_2(1+\mu||F(u_n)||_{\infty})||\nabla F(u_n)||_{\infty}||u_n||_{\infty}vol(M)
\endaligned
\]
Let
\[I(\lambda):=\sup\limits_{|z|\leqslant\lambda}|(HessF)(z)|,\]
\[J(\lambda):=\sup\limits_{|z|\leqslant\lambda}|\nabla F(z)|\]
and
\[H(\lambda):=\sup\limits_{|z|\leqslant\lambda}|F(z)|\]
they are both monotonously increasing functions and locally Lipschitz. So their derivatives exist almost everywhere and are not smaller than zero. We get
\[
\aligned
\frac{1}{2}\frac{d}{dt}\int_M|u_n|^2\,dM\leqslant&\gamma\cdot J(||u_n||_{\infty})||u_n||_{\infty}\sqrt{vol(M)}||\Delta u_n||_2\\
&+\kappa_2(1+\mu\cdot H(||u_n||_{\infty}))\cdot J(||u_n||_{\infty})||u_n||_{\infty}\cdot vol(M)
\endaligned
\]
using $\lambda^2_iC_i^n(t)$ to multiply two sides of (\ref{2.1}) and then summing $i$ from 1 to $n$, we get that:
\[
\aligned
&\frac{1}{2}\frac{d}{dt}\int_M|\Delta u_n|^2\,dM+\kappa_1\int_M|\nabla\Delta u_n|^2\,dM\\
=&\gamma\int_M(\nabla F(u_n)\times\Delta u_n)\cdot\Delta^2u_n\,dM-\kappa_2\int_M(1+\mu\cdot F(u_n))\nabla F(u_n)\cdot\Delta^2u_n\,dM\\
=&-\gamma\int_M[\nabla(\nabla F(u_n))\times\Delta u_n]\nabla\Delta u_n\,dM+\kappa_2\mu\int_M(\nabla F(u_n)\cdot\nabla u_n)(\nabla F(u_n)\cdot\nabla\Delta u_n)\,dM\\
&+\kappa_2\int_M(1+\mu\cdot F(u_n))\nabla(\nabla F(u_n))\cdot\nabla\Delta u_n\,dM\\
=&-\gamma\int_M\{[(HessF)(u_n)\cdot\nabla u_n]\times\Delta u_n\}\cdot\nabla\Delta u_n\,dM+\kappa_2\mu\int_M(\nabla F(u_n)\cdot\nabla u_n)(\nabla F(u_n)\cdot\nabla\Delta u_n)\,dM\\
&+\kappa_2\int_M(1+\mu\cdot F(u_n))[(Hess F)(u_n)\cdot\nabla u_n]\cdot\nabla\Delta u_n\,dM\\
\leqslant&\gamma\cdot I(||u_n||_{\infty})\int_M|\nabla u_n|\cdot|\Delta u_n|\cdot|\nabla\Delta u_n|\,dM+\kappa_2\mu\cdot J(||u_n||_{\infty})^2\int_M|\nabla u_n|\cdot|\nabla\Delta u_n|\,dM\\
&+\kappa_2(1+\mu\cdot H(||u_n||_{\infty}))\cdot I(||u_n||_{\infty})\int_M|\nabla u_n|\cdot|\nabla\Delta u_n|\,dM\\
\leqslant&\gamma\cdot I(||u_n||_{\infty})\cdot||\nabla u_n||_6\cdot||D^2u_n||_3\cdot||\nabla\Delta u_n||_2+\kappa_2\mu\cdot J(||u_n||_{\infty})^2\cdot C_1\cdot||\nabla u_n||_6\cdot||\nabla\Delta u_n||_2\\
&+\kappa_2(1+\mu\cdot H(||u_n||_{\infty}))\cdot I(||u_n||_{\infty})\cdot C_1\cdot||\nabla u_n||_6\cdot||\nabla\Delta u_n||_2
\endaligned
\]
By h\"{o}lder inequality:
\[
\aligned
&\frac{1}{2}\frac{d}{dt}(||\Delta u_n||_2^2)+\frac{\kappa_1}{2}||\nabla\Delta u_n||_2^2\\
\leqslant& C_2\cdot I(||u_n||_{\infty})^2\cdot||\nabla u_n||^2_6\cdot||D^2u_n||_3^2+C_2\cdot J(||u_n||_{\infty})^4||\nabla u_n||_6^2\\
&+C_2\cdot(1+\mu\cdot H(||u_n||_{\infty}))^2I(||u_n||_{\infty})^2||\nabla u_n||_6^2
\endaligned
\]
Due to $(2.4),(2.5),(2.6),(2.7)$ of \cite{CF}, there exist two locally lipschitz and monotonously increasing functions $Q(\lambda)$, $V(\lambda)$ such that
\[\frac{d}{dt}(||\Delta u_n||^2_2)\leqslant Q(||u_n||_2^2+||\Delta u_n||_2^2)\]
\[\frac{d}{dt}(||u_n||_2^2)\leqslant V(||u_n||_2^2+||\Delta u_n||_2^2)\]
So we get
\[\frac{d}{dt}(||\Delta u_n||^2_2+||u_n||_2^2)\leqslant(Q+V)(||\Delta u_n||^2_2+||u_n||_2^2)\]
Note that
\[||u_n(0)||_2\leqslant||u_0||_2\]
\[||\Delta u_n(0)||_2\leqslant||\Delta u_0||_2\]
Because of Exercise 3 in page 29 of \cite{T}, recalling that $Q+V$ is increasing and locally lipschitz, we get that there exists a $T^*$ and a $C_3(T^*)$ such that for all $t\in[0,T^*]$
\[||u_n(t)||_2^2+||\Delta u_n(t)||_2^2\leqslant C_3(T^*)\]
So by $(2.4)$ of \cite{CF}, we have that
\[||u_n(t)||_{\infty}\leqslant C\cdot(||u_n(t)||_2^2+||\Delta u_n(t)||_2^2)^{\frac{1}{2}}\leqslant C\sqrt{C_3(T^*)}\]
Multiplying two sides of (\ref{2.1}) by $\frac{d}{dt}C^n_i$ and summing $i$ from 1 to $n$, we get
\[
\aligned
||\partial_tu_n||^2_2=&\kappa_1\int_M\Delta u_n\cdot\partial_tu_n\,dM+\gamma\int_M\nabla F(u_n)\times\Delta u_n\cdot\partial_tu_n\,dM\\
&-\kappa_2\int_M(1+\mu\cdot F(u_n))\nabla F(u_n)\cdot\partial_tu_n\,dM\\
\leqslant&\frac{1}{4}||\partial_tu_n||_2^2+C_4\cdot||\Delta u_n||_2^2+\gamma\cdot J(||u_n||_{\infty})\cdot||\Delta u_n||_2\cdot||\partial_tu_n||_2\\
&+\kappa_2(1+\mu\cdot H(||u_n||_{\infty}))\cdot J(||u_n||_{\infty})\cdot||\partial_tu_n||_1
\endaligned
\]
Using H\"{o}lder inequality, we get
\[||\partial_tu_n||_2\leqslant C_5(T^*)\]
By $(2.1)$ of \cite{CF}, we obtain that
\[||u_n||_{W^{2,2}}\leqslant C\cdot(||u_n||_2^2+||\Delta u_n||^2_2)^{\frac{1}{2}}\leqslant C\sqrt{C_3(T^*)}\]
In conclusion:\\
$\bullet\,\{u_n\}$ is a bounded sequence in $L^{\infty}([0,T^*],W^{2,2})$;\\
$\bullet\,\{\partial_tu_n\}$ is a bounded sequence in $L^{\infty}([0,T^*],L^2)$.\\
Since $W^{2,2}\hookrightarrow C_B$ compactly and $W^{2,2}\hookrightarrow W^{1,2}$ compactly, we have:\\
\uppercase\expandafter{\romannumeral1}. $u_n\longrightarrow u$ strongly in $L^{\infty}([0,T^*],C_B)$\\
\uppercase\expandafter{\romannumeral2}. $u_n\longrightarrow u$ strongly in $L^{\infty}([0,T^*],W^{1,2})$\\
\uppercase\expandafter{\romannumeral3}. $\partial_tu_n\longrightarrow\partial_tu$ weakly* in $L^{\infty}([0,T^*],L^2)$\\
\uppercase\expandafter{\romannumeral4}. $\Delta u_n\longrightarrow\Delta u$ weakly* in $L^{\infty}([0,T^*],L^2)$\\
So
\[\partial_tu=\kappa_1\Delta u+\gamma\nabla F(u)\times\Delta u-\kappa_2(1+\mu\cdot F(u))\nabla F(u)\]
By Boundary Trace Imbedding Theorem, we have
\[
\aligned
&\Big|\Big|\frac{\partial u}{\partial\nu}\Big|\Big|_{L^2(\partial M)}\\
=&\Big|\Big|\frac{\partial u}{\partial\nu}-\frac{\partial u_n}{\partial\nu}\Big|\Big|_{L^2(\partial M)}\\
\leqslant&C||\nabla u-\nabla u_n||_{L^2(\partial M)}\\
\leqslant&\tilde{C}||\nabla u-\nabla u_n||_{L^2(M)}
\endaligned
\]
So
\[\frac{\partial u}{\partial\nu}=0\,\,\,\,a.e\,\,\,\,\partial M\]
This completes the proof.\endproof

Zonglin Jia

{\small\it Academy of Mathematics and Systems Sciences, Chinese Academy of Sciences, Beijing 100080,  P.R. China.}

{\small\it Email: 756693084@qq.com}\\

\end{document}